\documentclass{article}

\usepackage{
 amsmath,
 amsxtra,
 amsthm,
 amssymb,
 etex,
 mathrsfs,
 stmaryrd
 }
\usepackage[all]{xy}
\usepackage{hyperref}
\usepackage[T2A,T1]{fontenc}
\usepackage[OT2,T1]{fontenc}
\newcommand{\Zhe}{\mbox{\usefont{T2A}{\rmdefault}{m}{n}\CYRZH}}

\newtheorem{theorem}{Theorem}[section]
\newtheorem{lemma}[theorem]{Lemma}

\newtheorem{proposition}[theorem]{Proposition}
\newtheorem{corollary}[theorem]{Corollary}

\theoremstyle{definition}
\newtheorem{defn}[theorem]{Definition}
\newtheorem{remark}[theorem]{Remark}

\newcommand{\bd}{\begin{defn}}
\newcommand{\ed}{\end{defn}}
\newcommand{\bl}{\begin{lemma}}
\newcommand{\el}{\end{lemma}}
\newcommand{\bp}{\begin{proposition}}
\newcommand{\ep}{\end{proposition}}
\newcommand{\bt}{\begin{theorem}}
\newcommand{\et}{\end{theorem}}
\newcommand{\bc}{\begin{corollary}}
\newcommand{\ec}{\end{corollary}}
\newcommand{\br}{\begin{remark}}
\newcommand{\er}{\end{remark}}
\newcommand{\ba}{\begin{array}}
\newcommand{\ea}{\end{array}}
\newcommand{\bpf}{\begin{proof}}
\newcommand{\epf}{\end{proof}}

\newcommand{\Z}{\mathbb{Z}}
\newcommand{\Q}{\mathbb{Q}}
\newcommand{\Zp}{\mathbb{Z}_{p}}
\newcommand{\Qp}{\mathbb{Q}_{p}}
\newcommand{\Op}{\mathcal{O}}
\newcommand{\Ep}{E[p^{\infty}]}
\newcommand{\Epp}{E[\p^\infty]}
\newcommand{\Ap}{A[p^{\infty}]}
\newcommand{\al}{\alpha}
\newcommand{\be}{\beta}
\newcommand{\Ga}{\Gamma}
\newcommand{\ga}{\gamma}

\newcommand{\La}{\Lambda}

 \DeclareMathOperator{\Gal}{Gal}
 \DeclareMathOperator{\rank}{rank}
\DeclareMathOperator{\corank}{corank}

\newcommand{\cyc}{\mathrm{cyc}}

\newcommand{\mM}{\mathcal{M}}

\newcommand{\p}{\mathfrak{p}}

\newcommand{\ot}{\otimes}
\newcommand{\ilim}{\displaystyle \mathop{\varinjlim}\limits}

\newcommand{\coker}{\mathrm{coker}\,}

\newcommand{\lra}{\longrightarrow}

\newcommand{\ps}[1]{\llbracket #1 \rrbracket}

  \DeclareFontFamily{U}{wncy}{}
  \DeclareFontShape{U}{wncy}{m}{n}{<->wncyr10}{}
  \DeclareSymbolFont{mcy}{U}{wncy}{m}{n}
  \DeclareMathSymbol{\sha}{\mathord}{mcy}{"58}

\numberwithin{equation}{section}

\begin{document}
\title{On the structure of fine Mordell-Weil groups over a $\Zp$-extension and its intermediate subextensions}
 \author{Meng Fai Lim\footnote{School of Mathematics and Statistics, Key Laboratory of Nonlinear Analysis and Applications (Ministry of Education), Central China Normal University, Wuhan, 430079, P.R.China. E-mail: \texttt{limmf@ccnu.edu.cn}} }
\date{}
\maketitle

\begin{abstract} \footnotesize
\noindent
In this paper, we investigate the structure of the fine Mordell-Weil groups over the intermediate subextensions of a given $\Zp$-extension $F_\infty$ of $F$.


\end{abstract}

\section{Introduction}
In this paper, we always let $p$ denote an odd prime. Let $F$ be a number field and $A$ an abelian variety over $F$. The ($p$-primary) fine Selmer
group of $A$ is a much studied object in Iwasawa theory (for instances, see
\cite{CS05, LimFineDoc, LimFinePreprint, Wu}) and plays a central role in both the statement and proof of the Iwasawa main conjecture (see \cite{K, Kob, PR00}). Analogous to the usual Selmer group, the fine Selmer group $R(A/F)$ fits into the following short exact sequence
\[ 0\lra \mM(A/F)\lra R(A/F) \lra \Zhe(A/F_n)\lra 0, \]
where $\mM(A/F)$ and $\Zhe(A/F)$ denote the fine Mordell-Weil group and fine Tate-Shafarevich
group, respectively, as defined by Wuthrich \cite{WuTS}. These groups can be thought as the ``fine''
counterpart of the usual Mordell-Weil group and Tate-Shafarevich group. In recent years, these fine objects have garnered significant research attention, as evidenced in the works \cite{GK, LeiZ, LimFineMWanti, LQW, Ra}.

The goal of the paper is study the structure of the fine Mordell-Weil groups over the intermediate subextensions of a given $\Zp$-extension $F_\infty$ of $F$. For each $n$, we let $F_n$ be the unique subextension of $F_\infty/F$ such that $|F_n:F|=p^n$. The fine Mordell-Weil group $\mM(A/F_n)$ has a natural $\Zp[\Gal(F_n/F)]$-module structure. Denote by $\Ga$ the Galois group $\Gal(F_\infty/F)$ and write  $\La$ for $\Zp\ps{\Ga}$. Via the canonical surjection $\La\twoheadrightarrow \Zp[\Gal(F_n/F)]$, we may view $\mM(A/F_n)$ as a $\La$-module structure. Recall that a homomorphism $M \lra N$ of finitely generated torsion $\La$-modules is said to be a pseudo-isomorphism if its kernel and cokernel are finite. We denote this relation by $M\sim N$. Finally, let $\Phi_j$ denote the $p^j$-th cyclotomic polynomial. We are now in a position to state our main result.

\bt[Theorem \ref{fine MWgrowth}]
Let $A$ be a $g$-dimensional abelian variety defined over a number field $F$ and $F_\infty/F$ a $\Zp$-extension. Suppose further that every prime of $F$ above $p$ does not split completely in $F_\infty/F$. Then there exists a sequence of integers $\{s_n\}$ such that the following pseudo-isomorphism holds:
 \[\mM(A/F_n)^\vee \sim \bigoplus_{j= 0}^n (\La/\Phi_j)^{\oplus s_j}.\]
Furthermore, there exists an integer $n_0$ such that for every $n\geq n_0$, the integers $s_n$ satisfy the following estimates
\[ \max\{0, e_n-g|F:\Q|\} \leq s_n \leq \max\{0, e_n-1\},\]
where
\[ e_n = \left\{
           \begin{array}{ll}
            \displaystyle\frac{\rank_{\Z}A(F_n) - \rank_{\Z}A(F_{n-1}) }{\phi(p^n)} , & \mbox{if $n\geq 1;$} \\
            \\
            \rank_{\Z}A(F) , & \mbox{if $n=0$.}
           \end{array}
         \right.
\]
\et

We like to make a remark on the notation in the statement of the theorem. The $n$ appearing the sequence $\{s_n\}$ is consistent with the $n$ in $F_n$, as in the $n$th-layer of the $\Zp$-extension $F_\infty/F$, and how they fit into the asserted pseudo-isomorphism in the following form \[\mM(A/F_0)^\vee \sim(\La/\Phi_0)^{\oplus s_0}, \mM(A/F_1)^\vee \sim(\La/\Phi_0)^{\oplus s_0}\oplus (\La/\Phi_1)^{\oplus s_1}\] and so on. In other words, roughly speaking, $\mM(A/F_{n-1})^\vee$ and $\mM(A/F_n)^\vee$ differs (up to pseudo-isomorphism) by a certain number of $\La/\Phi_n$, and this number is given by $s_n$.

 Results of such have previously been established for an elliptic curve over the cyclotomic extension of $\Q$ (see \cite{LeiZ}) and certain $\Zp$-extensions of an imaginary quadratic field (see \cite{LimFineMWanti}). A key improvement over these earlier works, in addition to extending the setting to abelian varieties and more general $\Zp$-extensions, is that our approach does not require the finiteness of the fine Tate-Shafarevich groups, a condition assumed in the aforementioned references.

A natural question is whether one can determine the precise values of $s_n$, rather than merely obtaining bounds. This turns out to be not easy; indeed, we can compute the exact values of
$s_n$ only for elliptic curves under certain conditions. We illustrate this by the next result.

\bp[Proposition \ref{fine MWgrowth ellpitic equality}]
Let $E$ be an elliptic curve defined over a number field $\Q$ and $F$ a Galois extension of $\Q$ such that $\Gal(F/\Q)\cong (\Z/2\Z)^{\oplus m}$. Let $K_i$ ($i=1,...,2^m-1$) be all the quadratic extensions of $\Q$ contained in $F$, and $E_i$ the corresponding quadratic twist of $E$. We set $E_0=E$. Let $F_n$ be the intermediate extension of $F^\cyc/F$ such that $|F_n:F| = p^n$. Then there exists a sequence of integers $\{s_n\}$ such that  the following pseudo-isomorphism
 \[\mM(E/F_n)^\vee \sim \bigoplus_{j= 0}^n (\La/\Phi_j)^{\oplus s_j}.\]
holds for every $n$. Furthermore, the integers $s_n$ satisfy the following equalities
\[  s_n = \sum_{i=0}^{2^m-1}\max\{0, e_{i,n}-1\},\]
where
\[ e_{i,n} = \left\{
           \begin{array}{ll}
            \displaystyle\frac{\rank_{\Z}E_i(\Q_n) - \rank_{\Z}E_i(\Q_{n-1}) }{\phi(p^n)} , & \mbox{if $n\geq 1;$} \\
            \\
            \rank_{\Z}E_i(\Q) , & \mbox{if $n=0$.}
           \end{array}
         \right.
\]
\ep

For the proof, we will require to employ certain properties of Weil restriction (in the sense of Weil \cite{Weil}) as developed in \cite{MRS}. In the case when the elliptic curve has complex multiplication, we also have analogous results on obtaining precise values over certain non-cyclotomic $\Zp$-extensions (see Propositions \ref{fine CM MWgrowth split} and \ref{fine CM MWgrowth inert}).

It is then natural to ask whether our results can shed light on the structure of the fine Mordell-Weil group over a $\Zp$-extension. We address this question in the final section (see Section \ref{final section}), and this is the only place where the finiteness of fine Tate-Shafarevich groups will have to come into play in order for our results on intermediate subextensions to be able to transfer to the $\Zp$-extension.

We now give a brief overview of the paper. In Section \ref{fine section}, we recall the definition of the fine Selmer groups, fine Mordell-Weil groups and fine Tate-Shafarevich groups.
Moving to Section \ref{MW section}, we review relevant results on the variation of Mordell-Weil groups in a $\Zp$-extension, with the discussion relying heavily on a result of Lee \cite{Lee}. Section \ref{main section} is devoted to proving our main results (Theorem \ref{fine MWgrowth}) on the structure of the fine Mordell-Weil groups in the intermediate subextensions of a $\Zp$-extension. Subsequently, in Section \ref{CM section}, we present our results on CM elliptic curves, where we can determine the precise structure of the fine Mordell-Weil groups. In Section \ref{final section}, which is the last section of the paper, we examine the implications of our result on the structure of the fine Mordell-Weil group over a $\Zp$-extension.

\subsection*{Acknowledgement}
We thank Debanjana Kundu and Antonio Lei for their comments and questions on the paper.
The author is partially supported
by the Open Research Fund of Hubei Key Laboratory of Mathematical Sciences (Central China Normal University) No. MPL2025ORG020.

\section{Fine objects over $\Zp$-extension} \label{fine section}

In this section, we first review the definition of our fine objects. Let $A$ be an abelian variety defined over a number field $F$. Let $S$ be a finite set of primes of $F$ containing the primes above $p$, the bad reduction primes of $A$ and the infinite primes. Denote by $F_S$ the maximal algebraic extension of $F$ that is unramified outside $S$. For every extension $\mathcal{L}$ of $F$ contained in $F_S$, we write $G_S(\mathcal{L})=\Gal(F_S/\mathcal{L})$, and let $S(\mathcal{L})$ denote the set of primes of $\mathcal{L}$ above $S$.

Let $L$ be a finite extension of $F$ contained in $F_S$. Following Coates, Sujatha and Wuthrich \cite{CS05, Wu}, the fine Selmer group $R(A/L)$ of $A$ over $L$ is defined by the exact sequence
\[0\lra R(A/L)\lra H^1(G_S(L),\Ap)\lra \bigoplus_{v\in S(L)} H^1(L_v, \Ap).\]

Next, we introduce the fine Mordell-Weil group and fine Tate-Shafarevich group of Wuthrich \cite{WuTS}. To begin, the fine Mordell-Weil group $\mathcal{M}(A/L)$ is defined by
\[ \mathcal{M}(A/L) = \ker\Big(A(L)\ot_{\Zp}\Qp/\Zp \lra \bigoplus_{v\in S(L), v|p} A(L_v)\ot_{\Zp}\Qp/\Zp \Big).\]
The fine Tate-Shafarevich group $\Zhe(A/L)$ is then defined to be
\[ \Zhe(A/L) = \coker\Big( \mM(A/L)\lra R(A/L)\Big).\]
Consequently, we have the following short exact sequence
\begin{equation} \label{eqn fine short exact} 0 \lra \mM(A/L) \lra R(A/L) \lra \Zhe(A/L)\lra 0,
\end{equation}
where it can be shown that $\Zhe(A/L)$ injects into $\sha(A/L)[p^\infty]$ (see discussion in \cite[Section 2]{WuTS} or \cite[Section 3]{LimFineDoc}).

For a given $\Zp$-extension $F_\infty$ of $F$, we denote by $\Ga$ the Galois group of the extension $F_\infty/F$ and let $\La$ denote the resulting Iwasawa algebra over $\Zp$. We in turn denote by $F_n$ the fixed field of $\Ga_n$, where $\Ga_n$ is the unique subgroup of $\Ga$ of index $p^n$. The fine Selmer group of $A$ over $F_\infty$ is then defined to be $R(A/F_\infty) = \ilim_n R(A/F_n)$ which comes naturally equipped with a $\La$-module. The $\La$-modules $\mM(A/F_\infty)$ and $\Zhe(A/F_\infty)$ are defined via similar limiting processes.

\section{Mordell-Weil growth in $\Zp$-extension of number fields and local fields} \label{MW section}

In this section, we recall certain results on the variation of Mordell-Weil groups in a $\Zp$-extension $F_\infty/F$. As before, we write $\Ga:= \Gal(F_\infty/F)$. In what follows, all intermediate objects $A(F_n)\ot\Qp/\Zp$ and their Pontryagin dual will be viewed as $\Ga$-modules  via the surjection $\Ga \twoheadrightarrow \Gal(F_n/F)$ without further mention. We begin with the following result.

\bp \label{MWgrowth}
Let $A$ be an abelian variety defined over a number field $F$ and $F_\infty/F$ a $\Zp$-extension. Define
\[ e_n = \left\{
           \begin{array}{ll}
            \displaystyle\frac{\rank_{\Z}A(F_n) - \rank_{\Z}A(F_{n-1}) }{\phi(p^n)} , & \mbox{if $n\geq 1;$} \\
            \\
            \rank_{\Z}A(F) , & \mbox{if $n=0$.}
           \end{array}
         \right.
\]
Then we have pseudo-isomorphism
\[ \big(A(F_n)\ot \Qp/\Zp \big)^\vee \sim \bigoplus_{j=0}^n(\La/\Phi_j)^{\oplus e_j}.\]
\ep

\bpf
For the verification of the proposition, we refer readers to \cite[Proposition 3.8]{LimFineMWanti}.  We emphasize that the proof rely heavily on the results of Lee \cite{Lee}. Furthermore, it is important to highlight that the finiteness of the Tate-Shafarevich group is not assumed in the statement or the proof of the proposition.
\epf

For our subsequent discussion, we also require the following corresponding local result. The proof of which can be found in \cite[Proposition 3.9]{LimFineMWanti}.

\bp \label{MWgrowthlocal}
Let $\mathcal{K}$ be a finite extension of $\Qp$ and $\mathcal{K}_\infty$ a $\Zp$-extension of $\mathcal{K}$.
Let $A$ be a $g$-dimensional abelian variety defined over $\mathcal{K}$.
 Then we have
\[ \big(A(\mathcal{K}_n)\ot \Qp/\Zp \big)^\vee \sim \bigoplus_{j=0}^n(\La/\Phi_j)^{\oplus g|\mathcal{K}:\Q_p|},\]
where $\mathcal{K}_n$ is the subextension of $\mathcal{K}_\infty/\mathcal{K}$ with $|\mathcal{K}_n:\mathcal{K}|=p^n$.
\ep

We now discuss a semi-local variant of the preceding result, which will be needed for the subsequent parts of the paper. Let $v$ be a prime of $F$ lying above $p$, and let $\Ga_v$ denote the decomposition group of $\Ga:=\Gal(F_\infty/F)$ at some fixed prime of $F_\infty$ above $v$. We assume throughout that $|\Ga:\Ga_v| = p^{n_v}$ with $n_v<\infty$. In other words, the prime $v$ does not split completely in $F_\infty/F$. Fix a topological generator $\ga$ of $\Ga$. We may then define $\ga_v:=\ga^{p^{n_v}}$, which is a topological generator of $\Ga_v$. To distinguish between the two cases, we write $\La(\Ga)$ (resp., $\La(\Ga_v)$) for the Iwasawa algebra of $\Ga$ (resp., $\Ga_v$). The cyclotomic polynomials over $\La(\Ga)$ and $\La(\Ga_v)$ are then denoted by $\Phi_j(\ga)$ and $\Phi_j(\ga_v)$ respectively. A straightforward verification shows that
\[ \Phi_j(\ga_v)=\left\{
                   \begin{array}{ll}
                     \Phi_{j+n_v}(\ga), & \hbox{if $j \geq 1$;} \\
                     \displaystyle \prod_{k=0}^{n_v}\Phi_{k}(\ga), & \hbox{if $j=0$.}
                   \end{array}
                 \right.
 \]
For $n> n_v$, we have
\[ \La(\Ga)/\Phi_n(\ga) = \Big(\bigoplus_{i=0}^{p^{n_v}-1}\ga^i\La(\Ga_v)\Big)/\Phi_{n-n_v}(\ga_v) = \bigoplus_{i=0}^{p^{n_v}-1}\ga^i\Big(\La(\Ga_v)/\Phi_{n-n_v}(\ga_v)\Big)\cong \bigoplus_{i=0}^{p^{n_v}-1}\La(\Ga_v)/\Phi_{n-n_v}(\ga_v). \]
When $n=n_v$, we have
\[ \La(\Ga)/\prod_{k=0}^{n_v}\Phi_{k}(\ga) = \Big(\bigoplus_{i=0}^{p^{n_v}-1}\ga^i\La(\Ga_v)\Big)/\Phi_{0}(\ga_v) = \bigoplus_{i=0}^{p^{n_v}-1}\ga^i\Big(\La(\Ga_v)/\Phi_{0}(\ga_v)\Big)\cong \bigoplus_{i=0}^{p^{n_v}-1}\La(\Ga_v)/\Phi_{0}(\ga_v). \]

\bp \label{MWgrowthsemilocal}
Under the notation as above, we have
\[ \bigoplus_{w_n\mid v}\big(A(F_{n, w_n})\ot \Qp/\Zp \big)^\vee \sim \bigoplus_{j=0}^n(\La(\Ga)/\Phi_j(\ga))^{\oplus g|F_{n,v}:\Q_p|}\]
for every $n\geq n_v$.
\ep

\bpf
By Proposition \ref{MWgrowthlocal}, one has
\[ \bigoplus_{w_n\mid v}\big(A(F_{n, w_n})\ot \Qp/\Zp \big)^\vee \sim \bigoplus_{w_n\mid v}\bigoplus_{j=0}^{n-n_v}(\La(\Ga_{v})/\Phi_j(\Ga_{v}))^{\oplus g|F_{n,v}:\Q_p|}.\]
By applying the formulas before this proposition, the latter can be seen to be isomorphic to
\[\left(\La(\Ga)/\prod_{k=0}^{n_v}\Phi_{k}(\ga)\times \bigoplus_{j=1}^{n-n_v} \La(\Ga)/\Phi_{j+n_v}(\ga)\right)^{\oplus g|F_{n,v}:\Q_p|}.\]
The conclusion now follows.
\epf

\section{Main result} \label{main section}

We can now in position to prove our main result.

\bt \label{fine MWgrowth}
Let $A$ be a $g$-dimensional abelian variety defined over a number field $F$ and $F_\infty/F$ a $\Zp$-extension. Suppose further that every prime of $F$ above $p$ does not split completely in $F_\infty/F$. Then there exists a sequence of integers $\{s_n\}$ such that  the following pseudo-isomorphism holds:
 \[\mM(A/F_n)^\vee\sim \bigoplus_{j= 0}^n (\La/\Phi_j)^{\oplus s_j}.\]
Furthermore, there exists an integer $n_0$ such that for every $n\geq n_0$, the integer $s_n$ satisfies the following estimates
\[ \max\{0, e_n-g|F:\Q|\} \leq s_n \leq \max\{0, e_n-1\},\]
where
\[ e_n = \left\{
           \begin{array}{ll}
            \displaystyle\frac{\rank_{\Z}A(F_n) - \rank_{\Z}A(F_{n-1}) }{\phi(p^n)} , & \mbox{if $n\geq 1;$} \\
            \\
            \rank_{\Z}A(F) , & \mbox{if $n=0$.}
           \end{array}
         \right.
\]
\et

\bpf
By definition, we have the following exact sequence
\begin{equation}\label{fineMWexactseq}
  0\lra \mM(A/F_n) \lra A(F_n)\ot\Qp/\Zp \lra \bigoplus_{v\mid p}\bigoplus_{w_n\mid v} A(F_{n, w_n})\ot \Qp/\Zp.
\end{equation}
Taking Pontryagin dual, we obtain the exact sequence
\[ \bigoplus_{v\mid p}\big(\bigoplus_{w_n\mid v} A(F_{n, w_n})\ot \Qp/\Zp\big)^\vee \lra \big(A(F_n)\ot\Qp/\Zp\big)^\vee \lra \mM(A/F_n)^\vee \lra 0.\]
Taking Proposition \ref{MWgrowth} into account, the rightmost surjection implies that we have
\[ \mM(A/F_n)^\vee\sim \bigoplus_{j= 0}^n (\La/\Phi_j)^{\oplus s_j}\]
for some integers $s_j$. For every prime $v$ of $F$ above $p$, we let $\Ga_v$ denote the decomposition group of $\Ga:=\Gal(F_\infty/F)$ at some fixed prime of $F_\infty$ above $v$. We then set $n_0:=\max\{n_v\mid v|p\}$, where $p^{n_v} = |\Ga: \Ga_v|$.
Now, let $n$ be an arbitrary integer $\geq n_0$. In view of Propositions \ref{MWgrowth} and \ref{MWgrowthsemilocal}, by comparing the number of summands $\La/\Phi_n$ in the exact sequence, we obtain the lower bound of the theorem. For the upper bound, one simply note that if $e_n>0$, then  a non-torsion point $P\in A(F_n)-A(F_{n-1})$ is sent to a non-torsion point in $A(F_{n,w})-A(F_{n-1,w})$ for some $w$.
\epf

\br
(a) The inspiration towards proving the preceding theorem comes from a result of Wuthrich \cite[Corollary 7.3]{WuTS}, where he proves the following estimate
 \[ \max\{0, \rank_\Z E(F)- |F:\Q|\} \leq \corank_{\Zp}\mM(E/F) \leq \max\{0, \rank_\Z E(F)-1\}\]
for an elliptic curve over a number field $F$. Our proof approach is essentially to consider an $\La$-adic analogue of this result, with the underlying idea drawn from a work of Lei \cite{LeiZ}.

(b) Results of such were established in \cite[Corollary 3.8]{LeiZ} and \cite{LimFineMWanti} (also see \cite{GK} for some results in this direction), but in a less general setting and under the additional assumption that the fine Tate-Shafarevich groups are finite. 

(c) One might ask whether we can say anything if a prime above $p$ is totally split in $F_\infty/F$. In this context, as can be seen from the proof, we still have the upper bound. However, we are not able to derive a lower bound, as this part of the proof relies on Proposition \ref{MWgrowthsemilocal} which in turn requires that the prime does not split completely. \er

One can obtain a slightly better estimate than what was stated above. Namely, we have the following.

\bp \label{fine MWgrowth estimates}
Let $A$ be an abelian variety defined over a number field $F$ and $F_\infty/F$ a $\Zp$-extension. Assume that $A$ is isogenous to a product of abelian varieties $\prod_{i=1}^r A_i$, where each $A_i$ has dimension $g_i$. Suppose further that every prime of $F$ above $p$ does not split completely in $F_\infty/F$. Then there exists a sequence of integers $\{s_n\}$ such that the following pseudo-isomorphism
 \[ \mM(A/F_n)^\vee \sim \bigoplus_{j= 0}^n (\La/\Phi_j)^{\oplus s_j}\]
holds for every $n$.
Furthermore, there exists an integer $n_0$ such that for every $n\geq n_0$, the integer $s_n$ satisfies the following estimates
\[ \sum_{i=1}^r\max\{0, e_{i,n}-g_i|F:\Q|\} \leq s_n \leq \sum_{i=1}^r\max\{0, e_{i,n}-1\},\]
where
\[ e_{i,n} = \left\{
           \begin{array}{ll}
            \displaystyle\frac{\rank_{\Z}A_i(F_n) - \rank_{\Z}A_i(F_{n-1}) }{\phi(p^n)} , & \mbox{if $n\geq 1;$} \\
            \\
            \rank_{\Z}A_i(F) , & \mbox{if $n=0$.}
           \end{array}
         \right.
\]
\ep

\bpf
Under the assumption of the proposition, we have a pseudo-isomorphism
\[ \mM(A/F_n) \sim \bigoplus_{i=1}^r\mM(A_i/F_n)\]
for every $n$. The conclusion of the proposition then follows from this and applying Theorem \ref{fine MWgrowth} to each $\mM(A_i/F_n)$.
\epf

In particular, the above proposition implies that, in estimating the number of summands of $\La/\Phi_j$, it suffices to treat the case of a simple abelian variety.

For the remainder of the section, we consider a base-changing context, where we can obtain a better estimate. Let $A$ be a simple abelian variety defined over a number field $F$ and $F_\infty$ a $\Zp$-extension of $F$. Suppose that $L$ is a finite Galois extension of $F$ such that $L\cap F_\infty =F$ (this is a very mild condition; one can always replace $F$ by a larger extension in $F_\infty$ to ensure this). Denote by $L_n$ ($0\leq n\leq \infty$) the compositum of $F_n$ and $L$.
Since $L\cap F_\infty= F$, we have that $L_n\cap F_\infty = F_n$ and so the Galois groups $\Gal(L_n/F_n)$ naturally identified with one another. We shall use $G$ to denote these Galois groups. Let $B:=\mathrm{Res}^L_F(A)$, where $\mathrm{Res}^L_F(A)$ is the Weil restriction of scalar of $A$ from $L$ to $F$ (in the sense of \cite[\S1.3]{Weil}). From the functoriality property of Weil restriction, we have $A(L_n) = B(F_n)$ for every $n$.

The group ring $\Q[G]$ is semisimple. Furthermore, there is a decomposition
\[\Q[G] = \bigoplus_\rho \Q[G]_\rho\]
of a direct sum of minimal two-sided
ideals indexed by the irreducible rational representations $\rho$ of $G$. Here $\Q[G]_\rho$ is the $\rho$-isotypic component of $\Q[G]$, i.e., the sum of all left ideals of $\Q[G]$ isomorphic to $\rho$. Following \cite{MRS}, we define the $\rho$-twist of $A$ by
\[A_\rho:=(\Q[G]_\rho\cap\Z[G])\ot A.\]
By \cite[Theorem 2.1(i)]{MRS}, this is an abelian variety with dimension $g(\dim_{\Q}\rho)$. With these notation, we can now state the following result.

\bp \label{base change}
Retain the above settings.  Suppose further that every prime of $F$ above $p$ does not split completely in $F_\infty/F$. Then there exists a sequence of integers $\{s_n\}$ such that we have
 \[\mM(A/L_n)^\vee \sim \bigoplus_{j= 0}^n (\La/\Phi_j)^{\oplus s_j}.\]
Furthermore, there exists an integer $n_0$ such that for every $n>n_0$, the integer $s_n$ satisfies the following estimates
\[ \sum_\rho\max\{0, e_{\rho,n}-g(\dim_{\Q}\rho)|F:\Q|\} \leq s_n \leq \sum_\rho\max\{0, e_{\rho,n}-1\},\]
where
\[ e_{\rho,n} = \left\{
           \begin{array}{ll}
            \displaystyle\frac{\rank_{\Z}A_\rho(F_n) - \rank_{\Z}A_\rho(F_{n-1}) }{\phi(p^n)} , & \mbox{if $n\geq 1;$} \\
            \\
            \rank_{\Z}A_\rho(F) , & \mbox{if $n=0$.}
           \end{array}
         \right.
\]
\ep

\bpf
Recall that $A(L_n) = B(F_n)$ by the functoriality of Weil restriction. On the other hand, it follows from \cite[Theorem 4.5]{MRS} that $B(F_n)$ is isogenous to $\bigoplus_\rho A_\rho(F_n)$. Consequently, we have a pseudo-isomorphism
\[ \mM(A/L_n) \sim \bigoplus_\rho \mM(A_\rho/F_n)\]
for every $n$. Applying Theorem \ref{fine MWgrowth} to each $\mM(A_\rho/F_n)$ and subsequently summing up the estimates obtained, we obtain the conclusion of the proposition.
\epf

We specialize to the context of an elliptic curve, where we can obtain precise values.

\bp \label{fine MWgrowth ellpitic equality}
Let $E$ be an elliptic curve defined over a number field $\Q$ and $F$ a Galois extension of $\Q$ such that $\Gal(F/\Q)\cong (\Z/2\Z)^{\oplus m}$. Let $K_i$ ($i=1,...,2^m-1$) be all the quadratic extensions of $\Q$ contained in $F$, and $E_i$ the corresponding quadratic twist of $E$. We set $E_0=E$. Let $F_n$ be the intermediate extension of $F^\cyc/F$ such that $|F_n:F| = p^n$. Then there exists a sequence of integers $\{s_n\}$ such that  the following pseudo-isomorphism holds:
 \[ \mM(E/F_n)^\vee \sim \bigoplus_{j= 0}^n (\La/\Phi_j)^{\oplus s_j}.\]
Furthermore, for every $n$, we have the following equalities
\[  s_n = \sum_{i=0}^{2^m-1}\max\{0, e_{i,n}-1\},\]
where
\[ e_{i,n} = \left\{
           \begin{array}{ll}
            \displaystyle\frac{\rank_{\Z}E_i(\Q_n) - \rank_{\Z}E_i(\Q_{n-1}) }{\phi(p^n)} , & \mbox{if $n\geq 1;$} \\
            \\
            \rank_{\Z}E_i(\Q) , & \mbox{if $n=0$.}
           \end{array}
         \right.
\]
\ep

\bpf
Since $\Q^\cyc/\Q$ is totally ramified at $p$ and $F/\Q$ is a $2$-extension, we see that $F^\cyc/F$ is totally ramified at every prime of $F$ above $p$. Therefore, we may take $n_0=0$ in Proposition \ref{base change}. Finally, note that in the context of the current proposition, the twist of $E$ is precisely the corresponding quadratic twist of $E$.
\epf

\br
When $F=\Q$, the preceding proposition was established by Lei \cite{LeiZ}.
\er

\section{CM elliptic curves} \label{CM section}

In this section, we present some cases, where we can obtain precise values on the structure of the fine Mordell-Weil groups as in Proposition \ref{fine MWgrowth ellpitic equality}.

Throughout this section, $K$ will always denote an imaginary quadratic field with ring of integers $\Op_K$. Let $E$ be an elliptic curve defined over $K$ with complex multiplication given by $\Op_K$. We shall always assume that our elliptic curve has good reduction at every primes above $p$.

We first consider the case that the prime $p$ splits completely in $K/\Q$. Denote by $\mathfrak{p}$ and $\overline{\mathfrak{p}}$ the primes of $K$ above $p$. Since the class group of $K$ is finite, there exists an integer $h$ such that $\mathfrak{p}^h= \pi\Op_K$ and $\overline{\mathfrak{p}}^h= \overline{\pi}\Op_K$ for some $\pi\in\Op_K$. We then set
\[ \Epp:=\cup_{n\geq 1}E[\pi^n]\quad \mbox{and}\quad E[\overline{\mathfrak{p}}^\infty]:=\cup_{n\geq 1}E[\overline{\pi}^n]  \]
which are cofree $\Zp$-modules of corank $1$. Moreover, there is a natural decomposition
\[\Ep = \Epp \oplus E[\overline{\mathfrak{p}}^\infty]\]
of $\Gal(\overline{\Q}/K)$-modules. 

For a finite extension $L$ of $K$, we define the fine $\p$-Mordell-Weil group by
\[ \mathcal{M}_\p(E/L) = \ker\left( E(L)\ot_{\Op_K} K_\p/\Op_{K,\p} \lra \bigoplus_{w|\p} E(L_w)\ot_{\Op_{K,\p}} K_\p/\Op_{K,\p}\right),\]
where $\Op_{K,\p}$ is the ring of integers of $K_\p$. Note that since we are assuming that $p$ splits in $K/\Q$, we have identifications $K_\p\cong \Qp$, $\Op_{K,\p}\cong \Zp$, $K_{\overline{\p}}\cong \Qp$ and $\Op_{K,\overline{\p}}\cong \Zp$.

The fine $\overline{\p}$-Mordell-Weil group $\mathcal{M}_{\overline{\p}}(E/L)$ is defined similarly. One can then check easily that $\mathcal{M}(E/L)= \mathcal{M}_\p(E/L)\oplus \mathcal{M}_{\overline{\p}}(E/L)$. Here the identification is taken over $\Zp$, where $\Zp$ acts on $\mathcal{M}_\p(E/L)$ through the isomorphism $\Op_{K,\p}\cong \Zp$, while $\Zp$ acts on $\mathcal{M}_{\overline{\p}}(E/L)$ via the isomorphism $\Op_{K,\overline{\p}}\cong \Zp$.

\bp \label{fine CM MWgrowth split}
Retain setup as above. Let $K_\infty/K$ be a $\Zp$-extension. Let $F$ be a Galois extension of $K$ such that $\Gal(F/K)\cong (\Z/2\Z)^{\oplus m}$. Let $L_i$ ($i=1,...,2^m-1$) be all the quadratic extensions of $K$ contained in $F$, and $E_i$ the corresponding quadratic twist of $E$. We set $E_0=E$. Let $F_n$ be the intermediate extension of $F_\infty/F$ such that $|F_n:F| = p^n$. Then there exists a sequence of integers $\{s_n\}$ such that  the following pseudo-isomorphism holds:
 \[ \mM(E/F_n)^\vee \sim \bigoplus_{j= 0}^n (\La/\Phi_j)^{\oplus 2s_j}.\]
 Furthermore, there exists an integer $n_0$ such that for every $n>n_0$, the integer $s_n$ satisfies the following estimates
\[ s_n = \sum_{i=1}^r\max\{0, f_{i,n}-1\},\]
where
\[ f_{i,n} = \left\{
           \begin{array}{ll}
            \displaystyle\frac{\rank_{\Op_K}E_i(K_n) - \rank_{\Op_K}E_i(K_{n-1}) }{\phi(p^n)} , & \mbox{if $n\geq 1;$} \\
            \\
            \rank_{\Op_K}E_i(K) , & \mbox{if $n=0$.}
           \end{array}
         \right.
\]
  \ep

\bpf
To begin with, we claim that for every $\Zp$-extension $K_\infty$ of $K$, the primes of $K$ above $p$ do not split completely in $K_\infty/K$. Plainly, this is true if both primes $\p$ and $\overline{\p}$ are ramified in $K_\infty/K$. Suppose that that if one of them, say $\overline{\p}$, is unramified over $K_\infty/K$. Then $K_\infty$ is precisely the $\Zp$-extension of $K$ unramified outside $\p$. In this situation, it follows from \cite[Chap. Proposition 1.9(iii)]{deS} that the prime $\overline{\p}$ does not split in $K_\infty/K$. This establishes our claim.

Since the quadratic twist $E_i$ of $E$ also has complex multiplication by $\Op_K$, by the argument in Proposition \ref{fine MWgrowth ellpitic equality}, it suffices to show the theorem for the elliptic curve $E$ over $K_\infty$. By applying Pontryagin dual to
\[ 0\lra \mathcal{M}_\p(E/K_n) \lra  E(K_n)\ot_{\Op_K} K_\p/\Op_{K,\p} \lra \bigoplus_{w|\p} E(K_{n,w})\ot_{\Op_{K,\p}} K_\p/\Op_{K,\p},\] and subsequently tensoring the resulting sequence with $F_\p=\Qp$, we obtain another exact sequence
\[ \bigoplus_{w\mid \p} \big(E(K_{n,w})\ot_{\Op_{K,\p}} K_\p/\Op_{K,\p}\big)^\vee \lra \big(E(K_n)\ot_{\Op_K} K_\p/\Op_{K,\p}\big)^\vee \lra \mM_\p(E/K_n)^\vee \lra 0.\]
By an $\Op_K$-analogue of Propositions \ref{MWgrowth} and \ref{MWgrowthsemilocal}, we have
\[\big(E(K_n)\ot_{\Op_K} K_\p/\Op_{K,\p}\big)^\vee\sim \bigoplus_{j=0}^n \big(\La/\Phi_j\big)^{\oplus f_n}\]
and
\[\bigoplus_{w\mid \p} \big(E(K_{n,w})\ot_{\Op_{K,\p}} K_\p/\Op_{K,\p}\big)^\vee \sim \bigoplus_{j=0}^n \La/\Phi_j\]
for large enough $n$. On the other hand, if $f_n>0$, then a non-torsion point $P\in E(K_n)-E(K_{n-1})$ is sent to a non-torsion point in $E(K_{n,w})-E(K_{n-1,w})$ for some $w$. Hence this yields the required conclusion of the proposition.
\epf

We now prove the analog of the preceding proposition in the context when the prime $p$ is inert in $K/\Q$. We shall write $\p$ for the unique prime of $K$ above $p$. In this context, for every Galois extension $L$ of $K$, both the fine Mordell-group $\mM(E/L)$ and Mordell-group $E(L)$ carry a natural $\Op_K$-module structure. Moreover, for every prime $w$ of $L$ above $\p$, the group $E(L_w)$ has a $\Op_{K,\p}$-module structure.

\bp \label{fine CM MWgrowth inert}
Suppose that $E$ is an elliptic curve defined over an imaginary quadratic field $K$ with complex multiplication given by $\Op_K$, and that the prime $p$ is inert in $K/\Q$. Let $K_\infty/K$ be a $\Zp$-extension. Let $F$ be a Galois extension of $K$ such that $\Gal(F/K)\cong (\Z/2\Z)^{\oplus m}$. Let $L_i$ ($i=1,...,2^m-1$) be all the quadratic extensions of $K$ contained in $F$, and $E_i$ the corresponding quadratic twist of $E$. We set $E_0=E$. Let $F_n$ be the intermediate extension of $F_\infty/F$ such that $|F_n:F| = p^n$. Then there exists a sequence of integers $\{s_n\}$ such that  one has
 \[\mM(E/F_n)^\vee \sim \bigoplus_{j= 0}^n (\Op_{K,\p}\ps{\Ga}/\Phi_j)^{\oplus 2s_j}.\]
 Furthermore, there exists an integer $n_0$ such that for every $n\geq n_0$, the integer $s_n$ satisfies the following estimates
\[ s_n = \sum_{i=1}^r\max\{0, f_{i,n}-1\},\]
where
\[ f_{i,n} = \left\{
           \begin{array}{ll}
            \displaystyle\frac{\rank_{\Op_K}E_i(K_n) - \rank_{\Op_K}E_i(K_{n-1}) }{\phi(p^n)} , & \mbox{if $n\geq 1;$} \\
            \\
            \rank_{\Op_K}E_i(K) , & \mbox{if $n=0$.}
           \end{array}
         \right.
\]
  \ep

\bpf
As before, it suffices to show the result over $K_\infty/K$. Since $p$ is inert in $K/\Q$, there is only one prime $\p$ of $K$ over $p$, and so it must ramify in $K_\infty/K$.

As seen in the proof of \cite[Theorem 4.2]{LimFineMWanti}, we have isomorphisms
\[ E(K_n)\ot_{\Zp} \Qp/\Zp \cong E(K_n)\ot_{\Op_{K,\p}} K_\p/\Op_{K,\p} \]
and
\[ E(K_{n,w})\ot_{\Zp} \Qp/\Zp \cong E(K_{n,w})\ot_{\Op_{K,\p}} K_\p/\Op_{K,\p} \]
for every prime $w$ of $K_n$ above $p$. Therefore, the defining sequence of the fine Mordell-Weil group can be rewritten as
\[ 0\lra \mM(E/K_n) \lra E(K_n)\ot_{\Op_{K,\p}} K_\p/\Op_{K,\p} \lra \bigoplus_{w|p}E(K_{n,w})\ot_{\Op_{K,\p}} K_\p/\Op_{K,\p}\]
which are viewed as $\Op_{K,\p}[\Gal(K_n/K)]$-modules. Furthermore, the quadratic extension $K_\p/\Q_p$ is unramified, and so the cyclotomic polynomial $\Phi_n$ remains irreducible over $\Op_{K,\p}\ps{\Ga}$. Hence, by applying the $\Op_{K,\p}$-analogue of the results in Section \ref{MW section} and following the argument of the preceding proposition, we can obtain the conclusion of the present proposition. \epf

\br
One can of course formulate and prove analogues of Proposition \ref{base change} in these CM context. We shall leave the details to the interested readers.
\er

\section{Structure of fine Mordell-Weil groups over $\Zp$-extensions} \label{final section}

In this final section, we say how our results have implications on the structure of the fine Mordell-Weil group over the full $\Zp$-extension $F_\infty$. Recall that by \cite[Proposition 2.1.2]{Lee}, there is a $\La$-map
\[ \big(A(F_\infty)\ot\Qp/\Zp\big)^\vee \lra \La^{\oplus r} \oplus \Big(\bigoplus_{n=1}^\infty (\La/\Phi_n)^{\oplus a_n} \Big) \]
with finite kernel and cokernel. The proof of this result relies on the
fact that $A(F_n)\ot\Qp/\Zp$ is divisible. However, this latter fact does not always hold for the fine Mordell group $\mM(A/F_n)$ (see \cite[Section 7]{WuTS}). Nevertheless, by the structure theory of finitely generated $\Zp\ps{\Ga}$-module, we still have a $\Zp\ps{\Ga}$-homomorphism
\begin{equation}\label{structure theory}
  \mM(A/F_\infty)^\vee\lra \La^u\oplus \Big( \bigoplus_{i=1}^r\La/\Phi_{c_i}^{b_i} \Big)\oplus N
\end{equation}
with finite kernel and cokernel, where $u$ is a non-negative integer, $c_i$ and $b_i$ are some positive integers (but we allow $r=0$), and $N$ is a torsion $\La$-module whose support does not contain any cyclotomic polynomials.

We can now state the main result of this section.

\bp \label{fineMWgeneral1}
Let $A$ be a $g$-dimensional abelian variety defined over $F$ and let $F_{\infty}$ be a $\Zp$-extension of $F$. Retain the notation in (\ref{structure theory}).
Suppose that $\Zhe(A/F_n)$ is finite for every $n$. Then there exists a $n_0$ such that whenever $n\geq n_0$, we have the estimates
\[ \max\{0, e_n-g|F:\Q|\} \leq u + \#\{c_i~|~ c_i = n\} \leq \max\{0, e_n-1\},\]
where
\[ e_n = \left\{
           \begin{array}{ll}
            \displaystyle\frac{\rank_{\Z}A(F_n) - \rank_{\Z}A(F_{n-1}) }{\phi(p^n)} , & \mbox{if $n\geq 1;$} \\
            \\
            \rank_{\Z}A(F) , & \mbox{if $n=0$.}
           \end{array}
         \right.
\]
\ep

\bpf
For each $n$, we have the following commutative diagram
\[    \entrymodifiers={!! <0pt, .8ex>+} \SelectTips{eu}{}\xymatrix{
    0 \ar[r] & \mM(A/F_n)  \ar[d]^{\al_n} \ar[r] &  R(A/F_n)
    \ar[d]^{\beta_n} \ar[r] & \Zhe(A/F_n) \ar[d]^{g_n}  \ar[r] &0\\
    0 \ar[r]^{} & \mM(A/F_\infty)^{\Ga_n} \ar[r]^{} & R(A/F_\infty)^{\Ga_n} \ar[r] &
   \Zhe(A/F_\infty)^{\Ga_n}   & }
\]
with exact rows. By \cite[Theorem 3.3]{LimFineDoc}, the map $\be_n$ has finite kernel and cokernel for every $n$. Combining this with the finiteness of $\Zhe(A/F_n)$, we see that the map $\al_n$ has finite kernel and cokernel for every $n$. Taking the latter into account, the conclusion of the proposition follows from Theorem \ref{fine MWgrowth}.
\epf

There is a conjecture that asserts that the dual fine Selmer group $R(A/F_\infty)^\vee$ is a torsion $\La$-module. To the best of the author's knowledge, this conjecture, stated in such generality, first appears in the work Perrin-Riou \cite{PR00}, and was later studied in greater depth by Wuthrich \cite{WuThe, Wu}. For further theoretical and numerical evidence in support of this conjecture, the reader is referred to \cite[Corollary 3.5 and Remark 5.2]{LimFineDoc} and \cite[Theorem 3.9 and Section 6]{LimFinePreprint}.

Assuming the torsionness conjecture holds, we obtain the estimates
\[ \max\{0, e_n-g|F:\Q|\} \leq \#\{c_i~|~ c_i = n\} \leq \max\{0, e_n-1\}\]
for large enough $n$.

\br
Of course, one can always work with a slightly weaker assumption, namely that $\mM(A/F_\infty)^\vee$ is torsion over $\La$. (For some discussion on this latter torsionness, we refer readers to \cite{LQW}.) However, we also note that if one is willing to assume the finiteness of $\Zhe(A/F_n)$, then it is shown in \cite{LimFineDoc} that $\Zhe(A/F_\infty)^\vee$ is torsion over $\La$. Therefore, under the assumption of the finiteness of $\Zhe(A/F_n)$, we see that $R(A/F_\infty)^\vee$ is torsion over $\La$ if and only if $\mM(A/F_\infty)^\vee$ is so.
\er

Finally, we can say more under the assumption that $A(F_\infty)$ is finitely generated.

\bp \label{fineMWgeneral}
Let $A$ be a $g$-dimensional abelian variety defined over $F$ and let $F_{\infty}$ be a $\Zp$-extension of $F$.
Suppose that the following statements are valid.
\begin{enumerate}
  \item[$(a)$] $\Zhe(A/F_n)$ is finite for every $n$.
  \item[$(b)$] $A(F_\infty)$ is finitely generated over $\Z$.
\end{enumerate}
Then we have
\[ \mM(A/F_\infty)^\vee\sim \bigoplus_{j\geq 0} (\La/\Phi_j)^{\oplus s_j}, \]
and there exists a $n_0$ such that whenever $n\geq n_0$, we have the estimates
\[ \max\{0, e_n-g|F:\Q|\} \leq s_n \leq \max\{0, e_n-1\},\]
where
\[ e_n = \left\{
           \begin{array}{ll}
            \displaystyle\frac{\rank_{\Z}A(F_n) - \rank_{\Z}A(F_{n-1}) }{\phi(p^n)} , & \mbox{if $n\geq 1;$} \\
            \\
            \rank_{\Z}A(F) , & \mbox{if $n=0$.}
           \end{array}
         \right.
\]
\ep

\bpf
Under the finite generation of $A(F_\infty)$, it then follows from \cite[Theorem 2.1.2]{Lee} that
\[ A(F_\infty)^\vee \sim \bigoplus_{n\geq 0}\big(\La/\Phi_n\big)^{\oplus a_n}\]
for some integers $a_n$. Since $\mM(A/F_\infty)^\vee$ is a quotient of $A(F_\infty)^\vee$, we also have
\[ \mM(A/F_\infty)^\vee \sim \bigoplus_{n\geq 0}\big(\La/\Phi_n\big)^{\oplus s_n}\]
for some integers $s_n$. The conclusion now follows from this and the preceding proposition.
\epf

We end with a few remarks.

\br
(a) Proposition \ref{fineMWgeneral} was first shown in \cite[Corollary 3.8]{LeiZ} (also see \cite[Proposition 3.6]{LimFineMWanti}).

(b) At a first viewing of Proposition \ref{fineMWgeneral}, it may appear that we have not invoked the torsionness of $R(A/F_\infty)^\vee$ in the proof. Although this seems the case, we should note that the assumptions of the said proposition already imply the torsionness of $R(A/F_\infty)^\vee$ (see \cite[Corollary 3.5]{LimFinePreprint}).

(c) Of course, one can derive more precise results regarding the structure of the fine Mordell-Weil group within the setting of Propositions \ref{fine MWgrowth ellpitic equality}, \ref{fine CM MWgrowth split} and \ref{fine CM MWgrowth inert} which we leave to the interested readers to state these explicitly.
\er

\br
It is natural to ask what conclusions can be drawn without assuming the finite generation of $A(F_\infty)$. In this case, we can still obtain structural results for a certain quotient of $\mM(A/F_\infty)^\vee$, which is denoted by $Y_f(A/F_\infty)$ (see \cite[Definition 3.4]{LimFineMWanti} for the definition of this module; also see the proof of \cite[Proposition 3.7]{LimFineMWanti} where it is shown that this module is a quotient of $\mM(A/F_\infty)^\vee$), by  appealing to \cite[Proposition 3.6]{LimFineMWanti}.
\er

\end{document}